\documentclass{amsart}
\usepackage{latexsym,amssymb,amsmath,amsthm,amscd,graphicx}
\usepackage{esint} % needed for comand \fint
\usepackage{color}
\usepackage{comment}
\numberwithin{equation}{section}
\newtheorem{theorem}{Theorem}[section]
\newtheorem{lemma}[theorem]{Lemma}
\newtheorem{corollary}[theorem]{Corollary}

\theoremstyle{definition}

\newtheorem{definition}[theorem]{Definition}

\theoremstyle{remark}

\begin{document}

\title[On the reflexivity of the spaces]{On the reflexivity of the spaces of variable integrability and summability}

%  Only \author and \address are required; other information is
%  optional. Remove any unused author tags.

%  author one information
\author{Arash Ghorbanalizadeh}
\address{
Department of Mathematics, Institute for Advanced Studies in Basic Sciences (IASBS), Zanjan 45137-66731, Iran}
%\curraddr{}
\email{ghorbanalizadeh@iasbs.ac.ir ; gurbanalizade@gmail.com}
\thanks{}
\author{Reza Roohi Seraji}
\address{
Department of Mathematics, Institute for Advanced Studies in Basic Sciences (IASBS), Zanjan 45137-66731, Iran}
%\curraddr{}
\email{rroohi@iasbs.ac.ir}
\thanks{}
\author{Yoshihiro Sawano}
\address{Department of Mathematical Science,
Chuo University University,
Kasuga,, Bunkyo-ku 112-8551, Tokyo, Japan}
%\curraddr{}
\email{yoshihiro-sawano@celery.ocn.ne.jp}
\thanks{}
%  author two information
\author{}
\address{}
\curraddr{}
\email{}
\thanks{}

%  \subjclass is required.
\subjclass[2010]{Primary 46E30, Secondary 42B35}

\keywords{variable mixed Lebesgue-sequence space, reflexivity, K\"{o}the duality}
\date{}

\dedicatory{}

%  "Communicated by" -- provide editor's name; required.
\commby{}

%  Abstract is required.

\begin{abstract}
In this paper, we show that under the condition $ 1<p_-, q_-, p_+, q_+<\infty$,
the space $\ell^{q(\cdot)} (L^{p(\cdot)})$ is reflexive.  In this way
we give an answer
 to open problem posed by H\"ast\"o in 2017 about the reflexivity
 of the variable mixed Lebesgue-sequence spaces $\ell^{q(\cdot)} (L^{p(\cdot)})$.
 What is important here is that the dual space of $\ell^{q(\cdot)} (L^{p(\cdot)})$ is specified.
As its direct corollary,
we show that the corresponding Besov space $B^{s(\cdot)}_{p(\cdot)q(\cdot)}$ is reflexive.
\end{abstract}

\maketitle

\section{Introduction}

Variable Lebesgue spaces, as their name implies,
generalize
classical Lebesgue spaces,
replacing the constant exponent $p$ with a variable
exponent function $p(\cdot)$.
Their roots
go back to the
work of Orlicz \cite{WO}. The resulting Banach function space $L^{p(\cdot)}$,
called the Lebesgue space with the variable exponent
$p(\cdot)$,
enjoys many
properties similar to Lebesgue spaces,
but it also differs from
Lebesgue spaces
in surprising and subtle
ways. For this reason,
the study of variable Lebesgue spaces
has an intrinsic interest, but
its applications are also crucial.
For example,
variable Lebesgue spaces are applied to
nonlinear elastic mechanics \cite{Zhikov2},
electrorheological fluids \cite{Ruzicka}
or image restoration \cite{LiLiPi}.
Detailed exposition of various properties of
Lebesgue and Sobolev spaces with variable exponent
can be found in recent books \cite{Cruz,Die} and survey papers \cite{di,sam}.
We refer to the references therein for more complete overview.

Recently, related spaces with variable exponents have been considered.
See \cite{Nakai} for Hardy spaces and Campanato spaces.
Among others
Besov spaces with variable exponents
are considered in \cite{AlHa2010}.
 The variable exponent Besov space $B^{s(\cdot)}_{p(\cdot)q(\cdot)}$
 is
 defined via the variable mixed Lebesgue-sequence space
 $\ell^{q(\cdot)} (L^{p(\cdot)})$.
 Some of its properties have been established in \cite{AlHa2010, GG, KeVy}.
 The current paper proves the reflexivity
 of the space $\ell^{q(\cdot)} (L^{p(\cdot)})$.
 In this way we give an affirmative answer
 to open problem posed by H\"ast\"o \cite{has}.

We recall the notion of variable exponent Lebesgue spaces.
By a variable exponent we shall mean a measurable function $p(\cdot): \mathbb{R}^n \rightarrow [c,\infty]$,
where $c>0$. We denote the set of such functions by $\mathcal{P}_0.$
The subset of variable exponents with range $[1,\infty]$
is denoted
by $\mathcal{P}$.
Write
$p_{+}:= \mathop{\rm ess \; sup}\limits p(\cdot)$
and
$p_{-}:= \mathop{\rm ess \; inf}\limits p(\cdot)$.
For $p(\cdot)\in \mathcal{P}_0$, we define
$\Omega_{\infty}:=\{x\in \mathbb{R}^n: p(x)=\infty\}$ and
$\Omega_0:=\mathbb{R}^n\setminus \Omega_{\infty}$.
The variable exponent Lebesgue space $L^{p(\cdot)}$
consists of all measurable functions $f$ for which there exists $\lambda >0$ such that the modular
\[
	 \varrho_{L^{p(\cdot)}}(\lambda^{-1}f)
:=\int_{\Omega_0}\left(\lambda^{-1}|f(x)|\right)^{p(x)}
dx
+\left\|\lambda^{-1}f\right\|_{L^{\infty}(\Omega_{\infty})}
\]
is finite. Given $p(\cdot)\in \mathcal{P}_0$ and
$f \in L^{p(\cdot)}$,
we define the following Luxemburg norm
\[
\|f\|_{L^{p(\cdot)}}:=\inf\{\lambda>0: \varrho_{L^{p(\cdot)}}(\lambda^{-1}f)\leq 1\}.
\]

Next,
we recall the notion of the variable mixed space
$\ell^{q(\cdot)} (L^{p(\cdot)})$
defined in \cite{AlHa2010}. For $p(\cdot), q(\cdot) \in \mathcal{P}_0$
and a sequence $f=(f_{\nu})_{\nu=1}^\infty$
in $L^{p(\cdot)}$, we define
its modular
by
\begin{align}\label{P1}
\varrho_{{\ell}^{q(\cdot)}(L^{p(\cdot)})}\Big(f\Big)&:= \sum_{\nu=1}^{\infty}\inf\Big\{\lambda_{\nu}>0 : \varrho_{p(\cdot)}\left(\lambda_{\nu}^{-\frac{1}{q(\cdot)}}|f_{\nu}|\right)\le 1\Big\},
\end{align}
where we use the convention $\lambda^{\frac{1}{\infty}} = 1$
for $\lambda>0$.
Also, in the case of $q^{+}<\infty,$ this modular can be written as
\begin{equation}\label{P1a}
\varrho_{{\ell}^{q(\cdot)}(L^{p(\cdot)})}\Big(f\Big)
=\sum_{\nu=1}^{\infty}\||f_{\nu}|^{q(\cdot)}\|_{L^{\frac{p(\cdot)}{q(\cdot)}}}.
\end{equation}

The variable mixed Lebesgue-sequence space $\ell^{q(\cdot)} (L^{p(\cdot)})$ is defined as follows:
\[
{\ell}^{q(\cdot)}(L^{p(\cdot)}):=\Big\{f=(f_\nu)_{\nu=1}^\infty :
\mbox{there exists } \lambda>0
\mbox{ such that }
\varrho_{{\ell}^{q(\cdot)}(L^{p(\cdot)})}\left(\lambda^{-1}f\right)<\infty \Big\}.
\]
If $p(\cdot), q(\cdot) \in \mathcal{P}_0$ satisfy $q_+<\infty$, then,
as in
\cite[Theorem 3.8]{AlHa2010}, the space $\ell^{q(\cdot)} (L^{p(\cdot)})$ is
a quasi-normed space, i.e.,
\begin{align*}
\|f\|_{{\ell}^{q(\cdot)}(L^{p(\cdot)})}&:= \inf\left\{\lambda>0\,:\,
\varrho_{{\ell}^{q(\cdot)}(L^{p(\cdot)})}\left(\lambda^{-1}f\right)\le 1\right\}
= \inf\Big\{\lambda>0\,:\,
\sum_{\nu=1}^{\infty}\left\|\Big|\frac{f_{\nu}}{\lambda}\Big|^{q(\cdot)} \right\|_{L^{\frac{p(\cdot)}{q(\cdot)}}}\le 1\Big\}
\end{align*}
is a quasi-norm on $\ell^{q(\cdot)} (L^{p(\cdot)})$.
According to \cite[Theorem 1]{KeVy} and \cite[Theorem 3.6]{AlHa2010},
$\|\cdot\|_{{\ell}^{q(\cdot)}(L^{p(\cdot)})}$ is a norm,
if $p(\cdot), q(\cdot) \in \mathcal{P}$ satisfy either
one of the following conditions:
\begin{enumerate}
\item
$1 \le q(\cdot) \le p(\cdot) \le \infty $,
\item
$p_-\ge 1$
and $q(\cdot)=q_+=q_- \ge 1$ is constant,
\item
$\frac{1}{p(\cdot)} + \frac{1}{q(\cdot)}\le 1 $.
\end{enumerate}
Here and below, inequalities for exponents
are understood in the sense of almost everywhere.

The goal of this paper is to answer
an
open conjecture
posed in \cite{has}
in the affirmative.
The proof of Theorem \ref{theo3.1} is given in ${\mathbb R}^n$ but this restriction
does not lose
any generality because a similar argument works for an arbitrary
measurable set.

\begin{theorem}\label{theo3.1}
If
$ p(\cdot), q(\cdot) \in {\mathcal P}$
satisfy
$1<p_-,q_-,p_+,q_+<\infty$ and
one of the conditions $(1)$, $(2)$ or $(3)$ above,
then
$\ell^{q(\cdot)} (L^{p(\cdot)})$
is reflexive.
\end{theorem}

A direct corollary of Theorem \ref{theo3.1}
is that the corresponding Besov space $B^{s(\cdot)}_{p(\cdot)q(\cdot)}$
is reflexive.
Fix a pair
$(\phi,\Phi) \in {\mathcal S}^2$,
where ${\mathcal S}$ denotes the Schwartz function space,
so that
\begin{equation}\label{eq:220201-1}
{\rm supp}({\mathcal F}\phi) \subset\left\{\frac12 \le |\xi| \le 2\right\}, \quad
{\rm supp}({\mathcal F}\Phi) \subset\left\{|\xi| \le 2\right\}, \quad
|{\mathcal F}\Phi|^2
+
\sum_{\nu=1}^\infty
|{\mathcal F}\phi_{\nu}|^2=1.
\end{equation}
Here,
${\mathcal F}$ stands for the Fourier transform
and
$\phi_\nu:=2^{\nu n}\phi(2^\nu\cdot)$.
Set
$\phi_0:=\Phi$.
 Following \cite{AlHa2010},
 the variable exponent Besov space $B^{s(\cdot)}_{p(\cdot)q(\cdot)}$
 collects all distributions $f \in \mathcal{S}' $ such that
\[
\|f\|_{B^{s(\cdot)}_{p(\cdot)q(\cdot)}}^{\phi}
:= \left\|(2^{\nu s(\cdot)} \phi_{\nu} \ast f)_{\nu=0}^\infty
\right\|_{\ell^{q(\cdot)} (L^{p(\cdot)})}< \infty.
\]
Note that the norm does not depend
on the pair
$(\phi,\Phi) \in {\mathcal S}^2$
satisfying (\ref{eq:220201-1}):
Any choice of such a pair
$(\phi,\Phi) \in {\mathcal S}^2$ results in an equivalent semi-norm
\cite[Theorem 5.5]{AlHa2010}.

The next corollary is an immediate consequence of Theorem \ref{theo3.1}.
\begin{corollary}\label{cor}
Let
$p(\cdot), q(\cdot) \in {\mathcal P}$ and
$s(\cdot) \in L^\infty$
satisfy
$1<p_-,q_-,p_+,q_+<\infty$
as well as
one of the conditions $(1)$, $(2)$ or $(3)$ above.
Assume in addition that there exist
$c$, $p_\infty$, $q_\infty>0$
such that
\[
|p(x)-p_\infty|
+|q(x)-q_\infty|\le\frac{c}{\log(e+|x|)}, \quad
\]
\[
|p(x)-p(y)|
+|q(x)-q(y)|
+|s(x)-s(y)|\le\frac{c}{\log\left(e+\frac{1}{|x-y|}\right)}
\]
for all $x,y \in {\mathbb R}^n$.
Then
$B^{s(\cdot)}_{p(\cdot)q(\cdot)}$ is reflexive.
\end{corollary}
Once we prove Theorem \ref{theo3.1},
we can prove Corollary \ref{cor} with ease.
If we let
\[
A:f \in B^{s(\cdot)}_{p(\cdot)q(\cdot)} \mapsto
(2^{\nu s(\cdot)}\phi_\nu*f)_{\nu=0}^\infty
\in \ell^{q(\cdot)} (L^{p(\cdot)})
\]
and
\[
B: f=(f_\nu)_{\nu=0}^\infty
\in \ell^{q(\cdot)} (L^{p(\cdot)})
\mapsto
\sum_{\nu=0}^\infty \phi_\nu*[2^{-\nu s(\cdot)}f_\nu]
\in B^{s(\cdot)}_{p(\cdot)q(\cdot)},
\]
then $A$ is bounded thanks to the definition
of the norm of $B^{s(\cdot)}_{p(\cdot)q(\cdot)}$.
We claim that $B$ is also bounded.
In fact, if $l \ge 3$,
then by the support condition of $\phi$ and $\Phi$, we have
\[
\phi_l*B f=\sum_{\nu=l-3}^{l+3} \phi_l*\phi_\nu*[2^{-\nu s(\cdot)}f_\nu].
\]
Therefore,
as in the proof of \cite[Theorem 5.5]{AlHa2010},
if we let $\eta_{l}:=2^{l n}(1+2^l|\cdot|)^{-2n-2}$ for  $l \in {\mathbb Z}$,
then
\[
|\phi_l*B f|\le c\sum_{\nu=l-3}^{l+3} \eta_{l}*[2^{-\nu s(\cdot)}|f_\nu|]
\]
for some constant $c>0$.
A similar argument works for $l=0,1,2$ to have
\[
|\phi_l*B f|\le c\sum_{\nu=\max(l-3,0)}^{l+3} \eta_{l}*[2^{-\nu s(\cdot)}|f_\nu|]
\]
for some constant $c>0$.
Hence, thanks to
\cite[Lemma 4.7]{AlHa2010},
\[
\|B f\|_{B_{p(\cdot)q(\cdot)}^{s(\cdot)}}
\le c
\sum_{\nu'=-3}^3
\|(f_{\max(0,\nu'+l)})_{l=0}^\infty\|_{\ell^{q(\cdot)}(L^{p(\cdot)})}
\le c
\|f\|_{\ell^{q(\cdot)}(L^{p(\cdot)})}.
\]
Since $B \circ A={\rm id}_{B^{s(\cdot)}_{p(\cdot)q(\cdot)}}$,
it follows from the reflexivity of
$ \ell^{q(\cdot)} (L^{p(\cdot)})$
that
$B^{s(\cdot)}_{p(\cdot)q(\cdot)}$ is reflexive.

The rest of this paper is devoted to the proof of Theorem \ref{theo3.1}.
For the proof of Theorem \ref{theo3.1},
we use the following notation:
Denote by $P_N$
the projection mapping onto $U_N$,
where
$U_N$
consists of all elements $(f_\nu)_{\nu=1}^\infty$ such that
$f_\nu=0$ if $\nu>N$.
Note that $P_N$ will act on various spaces
but that $P_N$ is a bounded linear operator in any case.
\section{Proof of Theorem \ref{theo3.1}}

\subsection{A dense subspace of $\ell^{q(\cdot)} (L^{p(\cdot)})$}

Recall that we defined
$q_-={\rm ess}\inf q(\cdot)$.
It is known in \cite[Theorem 6.1(i)]{AlHa2010}
that
$\ell^{q_-}(L^{p(\cdot)})$
is embedded
into
$\ell^{q(\cdot)} (L^{p(\cdot)})$
with the embedding constant $1$.
We check that this embedding is dense.
\begin{lemma}\label{lem3.0}
Suppose $p_+,q_+<\infty$.
Then
$\ell^{q_-}(L^{p(\cdot)})$
is densely embedded
into
$\ell^{q(\cdot)} (L^{p(\cdot)})$.
\end{lemma}

\begin{proof}
It suffices to show that
$\bigcup\limits_{N=1}^\infty (\ell^{q(\cdot)}(L^{p(\cdot)}) \cap U_N)$
is
dense
in $\ell^{q(\cdot)} (L^{p(\cdot)}) $.
Take
$f=(f_\nu)_{\nu=1}^\infty \in \ell^{q(\cdot)}(L^{p(\cdot)})$
with norm less than $1$.
Then
\[
\sum_{\nu=1}^{\infty}\inf\Big\{\lambda_{\nu}>0 \,:\,
\varrho_{p(\cdot)}\left(\lambda_{\nu}^{-\frac{1}{q(\cdot)}}f_{\nu}\right)\le 1\Big\} < 1.
\]
We let $(\lambda_\nu)_{\nu=1}^\infty$ be
a sequence of non-negative real numbers
such that
$\sum\limits_{\nu=1}^\infty \lambda_\nu<1$
and that
$\varrho_{p(\cdot)}\left(\lambda_{\nu}^{-\frac{1}{q(\cdot)}}f_{\nu}\right)\le 1$.
Choose $\varepsilon \in (0,1)$ arbitrarily and
take $N$ large enough so that
$\sum\limits_{\nu=N+1}^{\infty}\lambda_\nu<
\varepsilon.
$
If we let
$g:=P_N f$,
then
$\|f-g\|_{\ell^{q(\cdot)} (L^{p(\cdot)})}<\varepsilon$.
Thus, we obtain the desired result.
\end{proof}

We equip
${\mathbb N} \times {\mathbb R}^n$
with the product measure
$\mu=\delta \times dx$,
where $\delta$ is the counting measure.
We can view
$\ell^{q(\cdot)} (L^{p(\cdot)})$
as a function space over
${\mathbb N} \times {\mathbb R}^n$.

\begin{definition}\label{defi:220130}
 For a vector-valued sequence space $S(X)$
 over a reflexive Banach space $X,$ define its K\"othe dual with
respect to the dual pair $(X, X^*)$ (see
\cite{GKP}) as follows:
\[
S(X)':=\left\{\bar{\varphi}=(\varphi_\nu)_{\nu=1}^\infty \in {X^{*}}^{\mathbb{N}}: \mbox{for each}~~ \bar{f}
=(f_\nu)_{\nu=1}^\infty \in S(X), \sum_{\nu=1}^\infty |\varphi_\nu(f_\nu)|<\infty \right\}.
\]
The norm of $\bar{\varphi}=(\varphi_\nu)_{\nu=1}^\infty \in S(X)'$ is given by
\[
\|\bar{\varphi}\|_{S(X)'}:=
\sup\left\{
\sum_{\nu=1}^\infty |\varphi_\nu(f_\nu)|\,:\,(f_\nu)_{\nu=1}^\infty \in S(X),
\bar{f}
=
\|(f_\nu)_{\nu=1}^\infty\|_{S(X)} \le 1
\right\}.
\]
\end{definition}

Recall that $p_->1$.
Under the identification
$(L^{p(\cdot)})^* \approx L^{p'(\cdot)}$
\cite{Die},
where $p'(\cdot)$ denotes the conjugate exponent,
it follows from Definition \ref{defi:220130} that
$(\ell^{q(\cdot)} (L^{p(\cdot)}))'$
stands for
the K\"{o}the dual of
$\ell^{q(\cdot)} (L^{p(\cdot)})$,
that is, a collection
$g=(g_\nu)_{\nu=1}^\infty \subset L^{p'(\cdot)}$
such that
$\displaystyle
\sum_{\nu=1}^\infty f_\nu g_\nu \in L^1
$
for all
$(f_\nu)_{\nu=1}^\infty \in \ell^{q(\cdot)} (L^{p(\cdot)})$.
The norm of such
$g$
is given by
\[
\|g\|_{(\ell^{q(\cdot)} (L^{p(\cdot)}))'}:=
\sup\left\{\left\|\sum_{\nu=1}^\infty f_\nu g_\nu \right\|_{L^1}\right\},
\]
where
$(f_\nu)_{\nu=1}^\infty \in \ell^{q(\cdot)} (L^{p(\cdot)})$
moves over all elements with norm less than or equal to $1$.

We investigate the continuity of $P_N$
on
$(\ell^{q(\cdot)} (L^{p(\cdot)}))'$.
\begin{lemma}\label{lem:210810-4}
For all
$g \in (\ell^{q(\cdot)} (L^{p(\cdot)}))'$,
$\lim\limits_{N \to \infty}P_N g=g$
in $(\ell^{q(\cdot)} (L^{p(\cdot)}))'$.
\end{lemma}

\begin{proof}
We must show that
$
\lim\limits_{N \to \infty}\|g-
P_N g\|_{(\ell^{q(\cdot)} (L^{p(\cdot)}))'}=0.
$
Here and below, we write $h_N:=g-P_N g$ for $N \in {\mathbb N}$.
Assume to the contrary that
$
\lim\limits_{N \to \infty}\|h_N\|_{(\ell^{q(\cdot)} (L^{p(\cdot)}))'}>0
$
keeping in mind that
$(\|h_N\|_{(\ell^{q(\cdot)} (L^{p(\cdot)}))'})_{N=1}^\infty
$
is a decreasing sequence.

Keeping in mind that $q_->1$,
fix $M>0$ so that
\begin{equation}\label{eq:210810-1}
\frac{1}{q_-}+\frac{1}{M q_-}<1.
\end{equation}

If we take $L \in {\mathbb N}$ large enough, then
\begin{equation}\label{eq:220131-101}
\|h_L\|_{(\ell^{q(\cdot)} (L^{p(\cdot)}))'}
\le
2^{\frac{1}{2M q_-}}
\lim_{N \to \infty}\|h_N\|_{(\ell^{q(\cdot)} (L^{p(\cdot)}))'}.
\end{equation}
From the definition of the space
$(\ell^{q(\cdot)} (L^{p(\cdot)}))'$,
there exists
$f=(f_\nu)_{\nu=1}^\infty \in \ell^{q(\cdot)} (L^{p(\cdot)})$
with norm less than $1$ such that
\begin{equation}\label{eq:220131-102}
\left\|\sum_{\nu=L+1}^\infty f_\nu g_\nu \right\|_{L^1}>
2^{-\frac{1}{2M q_-}}
\|h_L\|_{(\ell^{q(\cdot)} (L^{p(\cdot)}))'}.
\end{equation}
By multiplying each $f_\nu$ by $e^{i\rho_\nu(\cdot)}$, where $\rho_\nu(\cdot)$
is a real-valued measurable function,
we may assume that $f_\nu g_\nu \ge 0$ for each $\nu \in {\mathbb N}$.
By the dominated convergence theorem,
there exists $L'>L$ such that
\begin{equation}\label{eq:220131-103}
\left\|\sum_{\nu=L+1}^{L'} f_\nu g_\nu \right\|_{L^1}>
2^{-\frac{1}{2M q_-}}
\|h_L\|_{(\ell^{q(\cdot)} (L^{p(\cdot)}))'}.
\end{equation}
Once again from
the definition of
$(\ell^{q(\cdot)} (L^{p(\cdot)}))'$
we see that
there exists
$f'=(f'_\nu)_{\nu=1}^\infty \in \ell^{q(\cdot)} (L^{p(\cdot)})$
with norm less than $1$ such that
\begin{equation}\label{eq:220131-104}
\left\|\sum_{\nu=L'+1}^\infty f'_\nu g_\nu \right\|_{L^1}>
2^{-\frac{1}{2M q_-}}
\|h_{L'}\|_{(\ell^{q(\cdot)} (L^{p(\cdot)}))'}.
\end{equation}
By multiplying each $f'_\nu$ by $e^{i\rho'_\nu(\cdot)}$, where $\rho'_\nu(\cdot)$
is a real-valued measurable function,
we may assume that $f'_\nu g_\nu \ge 0$ for each $\nu \in {\mathbb N}$.
We let
$
f=(f''_{\nu})_{\nu=1}^\infty:=f'-P_{L'}f'+P_{L'}f-P_{L}f
$.
Then, we claim
\begin{equation}\label{eq:220131-105}
\|f''\|_{\ell^{q(\cdot)} (L^{p(\cdot)})}
\le 2^{\frac{1}{q_-}}<2.
\end{equation}
In fact,
if we write out the norm fully, then we have
\begin{align*}
\varrho_{{\ell}^{q(\cdot)}(L^{p(\cdot)})}\Big(2^{-\frac{1}{q_-}}f''\Big)
&=
\sum_{\nu=L+1}^{L'}\inf\Big\{\lambda_{\nu}>0 : \varrho_{p(\cdot)}\left(2^{-\frac{1}{q_-}}\lambda_{\nu}^{-\frac{1}{q(\cdot)}}f_{\nu}\right)\le 1\Big\}\\
&\quad+
\sum_{\nu=L'+1}^{\infty}\inf\Big\{\lambda_{\nu}>0 : \varrho_{p(\cdot)}\left(2^{-\frac{1}{q_-}}\lambda_{\nu}^{-\frac{1}{q(\cdot)}}f'_{\nu}\right)\le 1\Big\}.
\end{align*}
From the definition of $q_-$
and the fact that
$f$ and $f'$ have norm less than $1$, we have
\begin{align*}
\varrho_{{\ell}^{q(\cdot)}(L^{p(\cdot)})}\Big(2^{-\frac{1}{q_-}}f''\Big)
&\le
\sum_{\nu=L+1}^{L'}\inf\Big\{\lambda_{\nu}>0 : \varrho_{p(\cdot)}\left(2^{-\frac{1}{q(\cdot)}}\lambda_{\nu}^{-\frac{1}{q(\cdot)}}f_{\nu}\right)\le 1\Big\}\\
&\quad+
\sum_{\nu=L'+1}^{\infty}\inf\Big\{\lambda_{\nu}>0 : \varrho_{p(\cdot)}\left(2^{-\frac{1}{q(\cdot)}}\lambda_{\nu}^{-\frac{1}{q(\cdot)}}f'_{\nu}\right)\le 1\Big\}\\
&=
\frac12
\sum_{\nu=L+1}^{L'}\inf\Big\{\lambda_{\nu}>0 : \varrho_{p(\cdot)}\left(\lambda_{\nu}^{-\frac{1}{q(\cdot)}}f_{\nu}\right)\le 1\Big\}\\
&\quad+
\frac12
\sum_{\nu=L'+1}^{\infty}\inf\Big\{\lambda_{\nu}>0 : \varrho_{p(\cdot)}\left(\lambda_{\nu}^{-\frac{1}{q(\cdot)}}f'_{\nu}\right)\le 1\Big\}
\le 1.
\end{align*}
A direct consequence of
(\ref{eq:220131-105}) is that
\begin{equation}\label{eq:210810-2}
\left\|\sum_{\nu=L+1}^\infty f''_\nu g_\nu \right\|_{L^1}
\le
2^{\frac{1}{q_-}}\|h_L\|_{(\ell^{q(\cdot)} (L^{p(\cdot)}))'}.
\end{equation}
Meanwhile,
since
$\displaystyle
\left\|\sum_{\nu=L+1}^\infty f''_\nu g_\nu \right\|_{L^1}
=
\left\|\sum_{\nu=L+1}^{L'} f_\nu g_\nu \right\|_{L^1}+
\left\|\sum_{\nu=L'+1}^\infty f'_\nu g_\nu \right\|_{L^1}
$
thanks to the fact that $f''_\nu g_\nu \ge 0$ for each $\nu \in {\mathbb N}$, we
deduce
from
(\ref{eq:220131-101})--(\ref{eq:220131-104}) that
\begin{align}
\left\|\sum_{\nu=L+1}^\infty f''_\nu g_\nu \right\|_{L^1}
&>
2^{1-\frac{1}{2M q_-}}
\|h_{L'}\|_{(\ell^{q(\cdot)} (L^{p(\cdot)}))'}
\nonumber\\
\label{eq:210810-3}
&\ge
2^{1-\frac{1}{2M q_-}}
\lim_{N \to \infty}\|h_N\|_{(\ell^{q(\cdot)} (L^{p(\cdot)}))'}\\
\nonumber
&\ge
2^{1-\frac{1}{M q_-}}\|h_L\|_{(\ell^{q(\cdot)} (L^{p(\cdot)}))'}>0.
\end{align}
In view of (\ref{eq:210810-1}),
(\ref{eq:210810-2}) and (\ref{eq:210810-3}) contradict.
\end{proof}

We specify the dual space
$(\ell^{q(\cdot)} (L^{p(\cdot)}))^*$
 of
$\ell^{q(\cdot)} (L^{p(\cdot)})$
keeping in mind that
$(\ell^{q(\cdot)} (L^{p(\cdot)}))'$
is canonically embedded into
$(\ell^{q(\cdot)} (L^{p(\cdot)}))^*$.
\begin{lemma}
The dual space
$(\ell^{q(\cdot)} (L^{p(\cdot)}))^*$
is
canonically isomorphic to
$(\ell^{q(\cdot)} (L^{p(\cdot)}))'$.
\end{lemma}
\begin{proof}
%The dual space
%$(\ell^{q(\cdot)} (L^{p(\cdot)}))^*$
%%contains
%$(\ell^{q(\cdot)} (L^{p(\cdot)}))'$
%with the embedding norm $1$.
%So,
We must show that
$(\ell^{q(\cdot)} (L^{p(\cdot)}))^*$
is contained in
$(\ell^{q(\cdot)} (L^{p(\cdot)}))'$.

Let
$\ell \in (\ell^{q(\cdot)} (L^{p(\cdot)}))^*$.
Then
$\ell(P_N f) \to \ell(f)$
for all $f \in \ell^{q(\cdot)} (L^{p(\cdot)})$
thanks to Lemma \ref{lem3.0}.
Fix $N \in {\mathbb N}$ for the time being.
Since the space $\ell^{q(\cdot)}(L^{p(\cdot)}) \cap U_N$
is isomorphic to
$(L^{p(\cdot)})^N$,
according to \cite{Die},
we can uniquely find
$(g_\nu^N)_{\nu=1}^N \subset (L^{p'(\cdot)})^N$ such that
$\displaystyle
\ell(P_N(f))=
\int_{{\mathbb R}^n}
\sum_{\nu=1}^N g_\nu^N(x)f_\nu(x)dx
$
for all $f=(f_\nu)_{\nu=1}^\infty \in \ell^{q(\cdot)} (L^{p(\cdot)})$.
By the uniqueness of $g_\nu^N$,
we have
$g_\nu^N=g_\nu^{N'}$ for all
$1 \le \nu \le N \le N'$.
Thus, letting
$g:=(g_\nu^\nu)_{\nu=1}^\infty$,
we obtain an element
$g$ such that
\[
\sup\left\{\left|
\int_{{\mathbb R}^n}
\sum_{\nu=1}^N g_\nu^N(x)f_\nu(x)dx\right|
\,:\,f=(f_\nu)_{\nu=1}^\infty \in \ell^{q(\cdot)}(L^{p(\cdot)}) \cap U_N,
\|f\|_{\ell^{q(\cdot)}(L^{p(\cdot)})}=1\right\}
<\infty.
\]
Thus,
it follows that $g \in (\ell^{q(\cdot)} (L^{p(\cdot)}))'$.
In view of Lemma \ref{lem3.0},
$\bigcup\limits_{N=1}^\infty (\ell^{q(\cdot)}(L^{p(\cdot)}) \cap U_N)$ is dense in $\ell^{q(\cdot)} (L^{p(\cdot)})$.
Thus, $\ell$ is generated by $g \in (\ell^{q(\cdot)} (L^{p(\cdot)}))'$.
\end{proof}

We conclude the proof of Theorem \ref{theo3.1}.
Let
$Q \in (\ell^{q(\cdot)} (L^{p(\cdot)}))^{**}$,
so that
$Q$ is a bounded linear mapping from
$(\ell^{q(\cdot)} (L^{p(\cdot)}))^{*}
\approx(\ell^{q(\cdot)} (L^{p(\cdot)}))'$
to ${\mathbb C}$.

Let $g=(g_\nu)_{\nu=1}^\infty \in (\ell^{q(\cdot)} (L^{p(\cdot)}))^* \approx (\ell^{q(\cdot)} (L^{p(\cdot)}))'$.
Since $P_Ng \to g$
in $(\ell^{q(\cdot)} (L^{p(\cdot)}))'$
as $N \to \infty$
thanks to Lemma \ref{lem:210810-4},
it follows that
$Q \circ P_N(f) \to Q(f)$
for all $f \in \ell^{q(\cdot)} (L^{p(\cdot)})$
as $N \to \infty$
in $(\ell^{q(\cdot)} (L^{p(\cdot)}))^{**}$.
Since
the dual of $(L^{p'(\cdot)})^N$ is isomorphic to $(L^{p(\cdot)})^N$,
it follows that there exists a sequence
$(f_\nu^N)_{\nu=1}^\infty \in \ell^{q(\cdot)} (L^{p(\cdot)}) \cap U_N$
such that
$\displaystyle
Q \circ P_N(g)=
\int_{{\mathbb R}^n}
\sum_{\nu=1}^N g_\nu^N(x)f_\nu^N(x)dx.
$
Therefore,
since $f_\nu^N=f_\nu^{N'}$ for all $1 \le \nu \le N \le N'$,
we have an element
$f:=(f_\nu^\nu)_{\nu=1}^\infty$
satisfying
$\displaystyle
Q \circ P_N(g)=
\int_{{\mathbb R}^n}
\sum_{\nu=1}^N g_\nu^N(x)f_\nu(x)dx.
$
Since
$f \in (\ell^{q(\cdot)} (L^{p(\cdot)}))''=\ell^{q(\cdot)} (L^{p(\cdot)})$
according to \cite[p. 13, Theorem 2.9]{BeSh-text-98},
we conclude that $Q$ is induced by $f$.

\subsection*{Acknowledgements} {The authors would like to thank Prof. P. Gorka for sharing with us open problems which was presented at the conference ``Nonstandard growth phenomena". Also, the authors would like to thank professor Gorka for his valuable comments and feedback at the early stages of this research.}
Yoshihiro Sawano was partially supported by Grand-in-Aid for Scientific Research (C), No.\,19K03546, for Japan Society for the Promotion of Science.

%  Bibliographies can be prepared with BibTeX using amsplain,
%  amsalpha, or (for "historical" overviews) natbib style.
\bibliographystyle{amsplain}
%  Insert the bibliography data here.

\end{document}